\documentclass[11pt,a4paper,twoside]{article}

\usepackage{amssymb,amsfonts,amsmath,amsthm,euscript,bm}
\usepackage[english]{babel}
\usepackage[latin1]{inputenc}
\usepackage{dsfont}
\usepackage{color}
\usepackage{geometry}
\usepackage{graphicx}
\usepackage{subfigure}
\usepackage{multirow}
\usepackage[bf,small]{caption}

\def\i{\mathrm i}

\def\e{\mathrm e}

\def\d{\mathrm d}

\def\1{ \mathds{1}}
\def\R{\mathbb{R}}

\def\N{\mathbb{N}}
\def\P{\mathbb{P}}
\def\E{\mathbb{E}}

\def\R{\mathbb{R}}

\def\Z{\mathbb{Z}}

\def\v{\mbox{\rm Var}}

\def\cov{\mbox{\rm Cov}}

\newtheorem{theo}{Theorem}
\newtheorem{lem}{Lemma}
\newtheorem{prop}{Proposition}

\newtheorem{cor}{Corollary}
\newtheorem{Ex}{Example}
\newtheorem{Rem}{Remark}

\newcommand {\nn}{\nonumber}

\newcommand {\noi}{\noindent}

\begin{document}

\author{P. Doukhan\footnote{UMR AGM8088, University Cergy-Pontoise.}, \
A. Jakubowski\footnote{Nicolaus Copernicus University, Torun.},  \
 S.R.C. Lopes\footnote{Federal University of Rio Grande de Sul, UFRGS.}   \ and \
D. Surgailis\footnote{Vilnius University.}}

\title{Discrete-time trawl processes with long memory }

 \maketitle

\begin{abstract}
We introduce a class of discrete time stationary  trawl processes
taking real or integer
values and
written as sums of past values of independent `seed' processes on shrinking intervals (`trawl heights'). Related trawl processes
in continuous time were studied in Barndorff-Nielsen (2011) and Barndorff-Nielsen et al. (2014),
however in our case 
the i.i.d. seed processes can be very general and need not be infinitely divisible.
In the case when the trawl height decays with the lag as $j^{-\alpha}$ for some $1< \alpha < 2 $,
the trawl process exhibits long memory and its covariance decays as $j^{1-\alpha}$.
We show that under general conditions on generic seed process,
the normalized partial  sums of such trawl process may tend either to
a fractional Brownian motion 
or to an $\alpha$-stable
L\'evy process.

\end{abstract}

\medskip

\noi {\it Keywords:} trawl process, integer and continuous-valued time series, long memory,
fractional  Brownian motion, L\'evy process.
\\
{\it AMS Classification subjects 2010}\\
60G22 Fractional processes, including fractional Brownian motion
\\
60G51 Processes with independent increments; Lévy processes
\\
60G99 Trawl process
\\
60K99 Long range memory process

\section{Introduction}\label{sec1}

The present paper introduces a class of stationary random processes of  the form
\begin{equation}\label{g1}
X_k\ =\  \sum_{j=0}^\infty \gamma_{k-j}(a_j), \qquad k \in \Z
\end{equation}
where $\gamma_k = \{\gamma_k (u), u \in  \R\}$ are i.i.d. copies of a generic process  $\gamma = \{\gamma (u), u \in \R\}$ tending to zero
in probability as $u \to 0$,
and $a_j \in \R$ for $j \in \N, \, \lim_{j \to \infty} a_j =0$ are deterministic numbers.  Clearly,
\eqref{g1} includes the class of causal moving averages $X_k = \sum_{j=0}^\infty a_j \xi_{k-j} $ in i.i.d. r.v.s
$\{\xi, \xi_k\}$, which correspond to a trivial process $\gamma = \{\gamma (u) = \xi u, u \in \R\}$. In as follows, we call
$X = \{X_k, k \in \Z\}  $ the {\it trawl process} corresponding to the {\it seed process }   $\gamma = \{\gamma (u), u \in \R\}$ and
{\it trawl } $ a = \{a_j, j \ge 0\} $.  The above terminology is borrowed from Barndorff-Nielsen et al. (2014)
\cite{bar2014} which considered a related class of trawl processes in continuous time represented
as stochastic integrals
\begin{equation}\label{g11}
Y_t\ =\  \int_{(-\infty, t] \times \R} \1( x \in  (0, d_{t-s})) L(\d x, \d s), \qquad t \in \R
\end{equation}
where $L(\d x, \d s) $ is a homogeneous L\'evy measure on $\R^2 $, with independent values on disjoint sets, and
$\{d_t, t\in \R_+ \} $ is a deterministic function satisfying certain conditions. In the  case when this function
takes constant values $d_t = a_j,$ if $ t \in (j, j+1],$ for $ j = 0,1, \dots $, the discretized process $\{Y_k, k \in \Z\} $ in
\eqref{g11} coincides with $\{ X_k, k \in \Z\}$ in  \eqref{g1} with independent increment (L\'evy)  seed process 
$\Big\{\zeta (u) =  \int_{(0, u] \times (0,1]} L(\d x, \d s), u \in \R \Big\} $. Clearly, an integer-valued seed
process $\gamma = \{\gamma (u), u \in \R\}$
in \eqref{g1} results in an integer-valued trawl process  $\{ X_k, k \in \Z\}$, similarly as in the case
of continuous-time trawl processes of   \eqref{g11} studied in \cite{bar2014}. On the other hand,
the discrete-time set-up allows us to consider very general seed processes $\gamma$ which need not be infinitely divisible or
have independent increments as in  \cite{bar2014}.
\\
(\cite{bar2011}, page 22) note that trawl processes represent a flexible class of stochastic processes which can be used to
model serially dependent count data and other stationary time series, where the marginal
distribution and the autocorrelation structure can be modeled independently from each other.
Particularly, trawl processes can exhibit long memory or long-range dependence,
which is usually associated with the divergence of the covariance series: $\sum_{k\in \Z} |\cov  (X_0, X_k)| = \infty$,
see \cite{gir2012}, and which occurs in models  \eqref{g1} and \eqref{g11} when the trawl function
decays sufficiently slowly with the lag, see    \cite{bar2014} and \S\, \ref{sec2} below.    (\cite{bar2014}, figure 6) exhibit sample paths and autocorrelation graphs of integer-valued trawl process with long-memory trawl function showing a remarkably slow decay and a disagreement between true and sample
autocorrelations based on a very large sample length.
\\
The main question studied in this paper, which is also
one of the basic questions for statistical applications of trawl processes, is the rate of convergence and the limit distribution
of the sample mean. We prove that for trawl process with long-memory trawl function $a_j$ decaying as $j^{-\alpha}, 1< \alpha < 2$
this limit distribution is either $\alpha$-stable or Gaussian, moreover, a non-Gaussian stable limit is typical
for integer valued seed (and trawl) process, while a Gaussian limit occurs for `continuous' seed processes, e.g.
diffusions or stochastic volatility processes. We note that our non-Gaussian result contradicts
the conjecture in   (\cite{bar2014}, page 708) about a Gaussian partial sums limit for long-memory trawl process in  \eqref{g11}.
In particular, for a standard Poisson seed process $\gamma $ and $a_j \sim c_0 j^{-\alpha}, 1 < \alpha < 2 $
we obtain, with  $H={(3-\alpha)}/2$, a sequence of processes
$$
Z_n(t)=\frac 1{n^{H}}\sum_{j=1}^{[nt]}(X_j-\E X_j)
$$
whose second order moments converge  to
those of a fractional Brownian motion, $B_H$ with index $H$:
$$\lim_{n\to \infty}\cov (Z_n(s),Z_n(t))=  \cov (B_H(s),B_H(t)), \qquad \forall s,t.$$
Moreover, $Z_n(t)\to 0$ in probability (the process is  evanescent) but
$n^{H-\frac1\alpha}Z_n(t)$ converges to a non-trivial limit  which is an $\alpha$-stable L\'evy process.  (Note  $H-\frac1\alpha=\frac{(2-\alpha)(\alpha-1)}{2\alpha}>0$ since $1<\alpha<2$.)
\\
A similar phenomenon (convergence of the partial sums process to a L\'evy stable process)
occurs for a number of long-range dependent stationary processes with finite variance, see \cite{TaL1986},
\cite{TaWS1997},  \cite{Konst1998},   or \cite{miko2002}, \cite{Wi2003}, \cite{ls2003}, \cite{su2004},
\cite{kajt2008}, \cite{pils2014} and the references
therein, although in most of the literature  this convergence is limited to finite-dimensional distributions.
For $M/G/\infty$ queue  with heavy-tailed activity periods,
the adequate functional convergence was proved in 
\cite{ResnickVdB2000}.
Since the limiting  stable processes in these works have independent increments,
the above behavior is sometimes called `distributional short-range dependence' in contrast to `distributional long-range dependence'
occurring when the limit of the partial sums process has dependent increments.
See \cite{deh2002}, \cite{lps2005}.
See also
\cite{lif2014} for a nice discussion of stable and Gaussian limits under long-range dependence.


\section{Discrete-time trawl process}\label{sec2}

\subsection{Existence  of discrete-time trawl process}\label{sec21}

Let  $\gamma_k = \{\gamma_k (u), u \in  \R\}$ be i.i.d. copies of a generic {\it seed process}  $\gamma = \{\gamma (u), u \in \R\}$
with finite variance $g(u) = \v (\gamma (u))$ and mean $\mu(u) = \E \gamma (u)$
tending to zero as $u \to 0 $ so that $\gamma(0) = 0$ and $\gamma (u) \to_\P 0 $ as $u \to 0$.
A trawl $ a = \{a_j \ge 0, j \in \N\} $ is a deterministic sequence such that $\lim_{j \to \infty} a_j = 0$.
We shall assume that
 \begin{equation} \label{mug}
|\mu(u)| = 
{\cal O}(g(u)) \to 0 \quad  (u \to 0)
\end{equation}
and
\begin{equation}\label{ga}
\sum_{j=0}^\infty g(a_j) < \infty.
\end{equation}
The trawl process $X = \{ X_k, k \in \Z \} $   corresponding to trawl  $ a = \{a_j \ge 0, j \in \N\} $
and seed process $\gamma = \{\gamma (u), u \in \R\}$ is defined as
\begin{equation}\label{Xk}
X_k\ =\  \sum_{j=0}^\infty \gamma_{k-j}(a_j), \qquad k \in \Z.
\end{equation}
Let
\begin{equation}\label{rho1}
\rho(u,v) = \cov   (\gamma(u), \gamma(v)), \qquad (u,v \in \R )
\end{equation}
denote the covariance
function of the seed process $\gamma $. The following statement is an easy
consequence of the Kolmorogov three series theorem.

\begin{prop} Let conditions \eqref{mug} and \eqref{ga} be satisfied. Then
the series
in \eqref{Xk} converges a.s. and in mean square for any $k \in \Z $, and defines
a stationary process with mean $\E X_k =  \sum_{j=0}^\infty \mu(a_j)$ and covariance function
\begin{equation} \label{rhoX}
{\cov}  (X_0,X_k) \ = \ \sum_{j=0}^\infty \rho (a_j, a_{j+k}), \qquad k \in \N.
\end{equation}
\end{prop}
\noi
Clearly, if the seed process takes integer values: $\gamma(u) \in \Z, \, u \in \R$, this property also holds
for the trawl process: $X_k \in \Z \, \ (\forall \, k \in \Z) $.  The following examples
show that the class of trawl processes is very large.

\begin{Ex} [Random line seed process] \label{ex1} {\rm Let  $\gamma (u) = \xi u, u \in \R $, where $\xi $ is a r.v. with zero mean
and variance $\sigma^2<\infty $. Then $\mu(u) = 0,$  $g(u) = \sigma^2 u^2 $ and condition \eqref{ga} translates to
$\sum_{j=0}^\infty a^2_j < \infty $. Then
$X$ in \eqref{Xk} is a moving-average:
\begin{equation} \label{MA}
X_k = \sum_{j=0}^\infty a_j \xi_{k-j},
\end{equation}
where $\{\xi_k, k \in \Z\}$ are i.i.d. copies of $\xi$. }
\end{Ex}

\begin{Ex}  [Brownian motion seed process]  \label{ex2} {\rm Let $\gamma (u) = B(u), u \in \R_+ $, where $B$ is a Brownian motion with zero mean and covariance
$\E B(u) B(v) = u \wedge v $ and $a_j \ge 0$. Then $X$ in \eqref{Xk} is a stationary Gaussian process with zero  mean
and covariance
$\cov (X_0, X_k) = \sum_{j=0}^\infty a_j \wedge a_{k+j}, k \in \N$.
Particularly,
if $a_j = a^j, a \in (0,1)$ then 
$\cov   (X_0, X_k) = a^k/(1-a) $ and
$X$ in \eqref{Xk} agrees with an AR(1) process
written as a moving-average in \eqref{MA} with Gaussian innovations $\xi_k \sim {\cal N}(0,\sigma^2) $ and
$\sigma^2 = 1+ a $.  }

\end{Ex}

\begin{Ex} [Poisson and Bernoulli seed processes]  \label{ex3} {\rm Let $\gamma (u) = P(u)$,  $u \in \R_+ $, where $P$
is a Poisson process with mean $\mu(u) = u$, covariance
$\rho(u,v) = \cov  (P(u), P(v)) = u \wedge v $ and $a_j \ge 0$.  Then $X$ in \eqref{Xk} is a stationary process
with mean  $\E X_k = \sum_{j=0}^\infty a_j $ and the same covariance as in Example \ref{ex2}. Moreover,
$X_k $ takes integer values and has a Poisson marginal distribution with mean $\E X_0$.
\\
The above example can be generalized by considering a mixed Poisson seed process  $\gamma(u)=P(u\zeta)$, where $P$ is as above and
$\zeta > 0$ is a random variable with $\E \zeta < \infty $,
independent of $P$. Particularly, \cite{fok}
proved that when $\zeta$ is  exponentially distributed then   $P(u\zeta)$ has  negative binomial  marginal distribution. 
\\
The Bernoulli seed process is defined by $\gamma(u)={\1}(U \le u) $, where  $U\sim {\cal U}[0,1]$ is a uniformly distributed random variable.
Note also
\begin{align*}
\mu(u)&=u\E\zeta, &\rho(u,v)&= (u \wedge v)\E\zeta + uv \v(\zeta)\ &\mbox{($\gamma $ \ is a mixed Poisson process)}, \\
\mu(u)&=u, &\rho(u,v) &= u\wedge v-uv \  &\mbox{($\gamma $ \ is a Bernoulli process).}
\end{align*}
}
\end{Ex}

\noi Further examples of trawl processes can  be found in \S~\ref{31} (Examples \ref{Ex4}-\ref{Ex5}) and \S~\ref{32} (Example \ref{Rem1}).  As explained in \S~\ref{sec1}, this paper is focused on long memory properties
and the behavior of the partial sums process of stationary trawl process $X$ in \eqref{Xk}.

\subsection{Second order properties of discrete-time trawl process}\label{sec22}

 The covariance
function $\cov  (X_0,X_k)$
in \eqref{rhoX} depends both on the trawl
$a = \{ a_j \} $ and on the covariance function $\rho (u,v) $ of the seed process. In order to characterize
long memory property in terms of the trawl  $a = \{ a_j \} $ alone, it is convenient
to impose a linear growth condition on the variance $g (u) = \v    (\gamma(u))$ at the origin $u=0$:
\begin{equation}\label{g0}
g(u) = |u|( 1 + o(1)), \qquad  u \to 0.
\end{equation}
Under \eqref{g0},  condition \eqref{ga} is equivalent to the summability
of the trawl sequence:
\begin{equation}\label{G12}
\sum_{j=0}^\infty |a_j| < \infty.
\end{equation}
Moreover, for obtaining more precise decay of the covariance function in \eqref{rhoX} we also assume that
\begin{eqnarray} \label{G11}
\rho(u,v)&=&(|u| \wedge |v|) (1 + o(1)),  \quad \text{as} \quad u, v \to  0,  \  uv  >0.
\end{eqnarray}
Clearly, the trawl processes in Examples \ref{ex2} and \ref{ex3} satisfy \eqref{g0}  and \eqref{G11}
provided the seed processes in these examples are suitably extended to negative $u <0$. 
Denote by $S_n = \sum_{k=1}^n X_k$ the partial sums process of the trawl process in
\eqref{Xk}.

\begin{prop} \label{propLM}  
(i) Assume conditions  \eqref{mug}, \eqref{g0},  \eqref{G11} and
\begin{equation} \label{g4}
a_j \ = \  c_0 j^{-\alpha} (1+ o(1)), \quad j \to \infty  \quad (\exists \, c_0 \ne 0, \, 1 < \alpha < 2).
\end{equation}
Then
\begin{equation} \label{covX}
\cov   (X_0, X_k) \ = \  c_1 k^{1-\alpha}(1+ o(1)), \qquad k \to \infty
\end{equation}
and
\begin{equation}\label{g5}
\v  (S_n)  = \sum_{k, l=1}^n  \cov  (X_k, X_l) \ \sim \ c_2 \, n^{3-\alpha} \gg n, \quad n \to \infty,
\end{equation}
where  $c_1 =  c_0/(\alpha -1)$, and $ c_2 = 2c_1/(2-\alpha)(3-\alpha)$.

\medskip

\noi (ii) Assume conditions \eqref{mug}, \eqref{g0},
\begin{equation} \label{rhobdd}
|\rho(u,v)| \le C (|u|\wedge |v|) \qquad (u,  v \in \R)
\end{equation}
and
\begin{equation} \label{G4}
\sum_{j=1}^\infty  j |a_j| \ < \ \infty.
\end{equation}
Then
\begin{equation} \label{G5}
\sum_{k=1}^\infty |\cov   (X_0, X_k)| \ < \ \infty
\end{equation}
and
\begin{equation}\label{g6}
\v  (S_n)  = n\sum_{|k|<n}\Big(1-\Big|\frac kn\Big|\Big)  \cov  (X_k, X_0)  \sim  \sigma^2 \, n,
\end{equation}
where $\sigma^2  =
\sum_{k\in \Z} \cov   (X_0, X_k)$.

\end{prop}

\begin{Rem} {\rm
The estimation of the parameter of interest needs additional work: it  will be considered in further papers. }
\end{Rem}


\noindent {\it Proof.} (i) Let $c_0 >0$ in \eqref{g4}, the case $c_0 < 0$ follows analogously. Then
$a_j >0$, and $a_{k+j} >0$ hold for all $k \ge 1 $ and  $j > j_0$, where $j_0$ is large enough.
Moreover, for any $\epsilon >0$ there exists $j_0< j_\epsilon < \infty $ such that
\begin{equation} \label{ajk}
a_{j+k}  <  a_j,    \quad  \text{for all} \quad  \forall \ j_\epsilon <  j < k/2\epsilon, \quad \forall \, k \ge 2\epsilon j_\epsilon.
\end{equation}
Indeed, by  \eqref{g4} we have that for any $\epsilon >0$ there exists $j_\epsilon > j_0> 0$ such that
$a_j >  c_0j^{-\alpha}(1-\epsilon)$, $a_{k+j} < c_0 (j+k)^{-\alpha} (1+ \epsilon)$ and therefore
$$
\Big( \frac{a_{j+k}}{a_j}\Big)^{\frac1\alpha} <  \frac j{j+k}\Big(\frac{1+ \epsilon}{1-\epsilon}\Big)^{\frac1\alpha}, \qquad \forall \,
j > j_\epsilon, \quad \forall \, k\ge 1.
$$
Since $((1+ \epsilon)/(1-\epsilon))^{\frac1\alpha}  < 1 + 2\epsilon$ if
$\epsilon >0$ is small enough, relation \eqref{ajk} follows
since   $ j/(j+k) \le 1/(1 + 2\epsilon)$ for $1 \le j <  k/2\epsilon $. 
Next, for sufficiently large $k$ ($k > 2\epsilon j_\epsilon$) split
$k^{\alpha -1} \cov   (X_0, X_k)
= \sum_{j=0}^\infty k^{\alpha -1} \rho(a_j, a_{k+j}) = \sum_{i=1}^3 I_{i,k}$, where
$$
I_{1,k} = \sum_{0\le j \le j_\epsilon} \dots, \qquad
I_{2,k} = \sum_{j_\epsilon < j  < k/2\epsilon} \dots, \qquad
I_{3,k} = \sum_{j \ge k/2\epsilon } \dots.
$$
By  \eqref{g0}, \eqref{g4} and Cauchy-Schwartz inequality, for any fixed $\epsilon >0$ and
$1 \leq j \leq j_\epsilon$,
$$
|\rho(a_j, a_{k+j})| \le g(a_j)^{\frac12} g(a_{k+j})^{\frac12}
\le C |a_{k+j}|^{\frac12} \le C k^{-\frac\alpha2}, \qquad k \to \infty
$$
implying 
$$
|I_{1,k}| \le C k^{\alpha-1} k^{-\frac \alpha2} = O(k^{ -(1-\frac \alpha2)}) = o(1), \qquad k \to \infty.
$$
Next, by \eqref{G11} and \eqref{g4}, $|\rho(a_j, a_{j+k})| \le C |a_j| \wedge |a_{j+k}|
\le  C j^{-\alpha}, \, (\forall \, j, k \ge 1) $
and therefore
$$
I_{3,k} \le C k^{\alpha -1} \sum_{j \ge k/2\epsilon } j^{-\alpha}
\le C \epsilon^{\alpha -1 }
$$
can be made arbitrarily small uniformly in $k \ge 1 $ by choosing
$\epsilon >0$ small enough. Finally, by \eqref{ajk} and \eqref{G11},
\begin{equation} \label{Ik}
I_{2,k} = c_0 k^{\alpha-1} \sum_{j_\epsilon < j < k/2\epsilon } \frac{1 + \delta_{j,k}}{ (k+j)^{\alpha}} ,
\end{equation}
where $\sup_{j \ge 1} |\delta_{j,k}| = 0$ as $k \to \infty $. Note that for each
$\epsilon >0$, as $k \to \infty $
\begin{eqnarray} \label{Jk}
J_{k}(\epsilon)&:=&k^{\alpha-1} \sum_{j_\epsilon < j < k/ 2\epsilon} (k+j)^{-\alpha} \
= \ \frac{1}{k} \sum_{\frac {j_\epsilon}k < \frac jk < 1/2\epsilon} \frac{1}{ \left(1 + \frac jk\right)^\alpha}  \nn \\
&\to&\int_0^{1/2\epsilon } \frac{\d x}{ (1 + x)^\alpha}\ = \ \frac{1}{\alpha -1} \left(1 - (2\epsilon)^{\alpha -1} \right).
\end{eqnarray}
According to \eqref{Ik} and \eqref{Jk}, for any $\delta >0$ and any $\epsilon_0>0$
one can find
$0<\epsilon < \epsilon_0 $  and  $K_0  >0 $
such that $|I_{2,k} - c_0/(\alpha -1)| < \delta $ holds for all $k > K_0$. This proves \eqref{covX} while
\eqref{g5} follows from   \eqref{covX}, see e.g. (\cite{gir2012}, proposition 3.3.1).
\smallskip

\noi (ii) It suffices to prove \eqref{G5} since \eqref{g6} follows from  \eqref{G5} and the dominated
convergence theorem.
According to \eqref{rhoX}, \eqref{rhobdd}, \eqref{G4},
\begin{eqnarray*}
\sum_{k=1}^\infty |\cov  (X_0,X_k)|
&\le&C\sum_{k=1}^\infty \sum_{j=0}^\infty |a_j| \wedge |a_{j+k}| \\
&\le&C  \sum_{k=1}^\infty \sum_{j=0}^\infty |a_{j+k}|  \ \le \ C \sum_{k=1}^\infty  k |a_k| \ < \ \infty.
\end{eqnarray*}
Proposition \ref{propLM} is proved. \hfill $\Box$


\section{Partial sums limits of trawl processes}

We shall consider two typical cases of the seed  process $\gamma $ in \eqref{g1}:

\begin{description}

\item [{\bf Case 1}:] $\gamma(u), u \ge 0 $ is centered: $\mu (u) = 0$  and a.s. continuous (e.g., a Brownian motion).

\item[ {\bf Case 2}:] $\gamma(u), u \ge 0 $ is a pure jump  process
(a typical example is a Poisson process with $\mu (u) = g(u) = u $).

\end{description}

\noindent
Particularly, in Examples \ref{ex2} and \ref{ex3} of $\gamma $ (Brownian motion and Poisson process) and a regularly
decaying trawl $a = \{a_j\}$ in \eqref{G4} with exponent $1< \alpha < 2 $
the conditions of Proposition \ref{propLM}~(i) are satisfied and the covariance function
of the trawl process decays as $ k^{1-\alpha }$, see  \eqref{covX}. The last fact
implies that  the variance of  $S_n = \sum_{k=1}^n X_k  $ grows faster than $n$, see \eqref{g6}.

\noi
In the following subsections we detail conditions on the seed process $\{\gamma (u), u \in \R\}$ which guarantee
that the partial sums process of the trawl process $\{X_k\}$ with regularly decaying trawl \eqref{g4}
tends to either a Gaussian process (fractional Brownian motion with Hurst parameter $H = (3-\alpha)/2 \in (1/2, 1) $ (Case 1) or to a
$\alpha$-stable L\'evy process
(Case 2).
\

\noindent
The following decomposition of the partial sums process as a sum of   independent random variables is crucial
for the proofs of Theorem \ref{thmgauss} and 
Theorem \ref{thmstable}.
\begin{lem} [Decomposition] \label{decomplem}
We have
\begin{eqnarray}\label{Zn}
S_n&=&\sum_{k=1}^n X_k \ = \ \sum_{s=-\infty}^n Z_{s,n},  \quad \text{where} \quad
Z_{s,n}\ =\ \sum_{k=1\vee s}^n \gamma_s(a_{k-s}).
\end{eqnarray}
Then the random variables  $(Z_{s,n})_{s \le n}$ are independent.
\end{lem}
Write $\to_{f.d.d.}$ for the weak convergence of finite-dimensional distributions and
$\to_{{\mathcal D}(J_1)}$ and $\to_{{\mathcal D}(M_1)} $ for the weak convergence
of random elements in the Skorohod space $D[0,1]$ endowed with the $J_1$-topology and the $M_1$-topology,
respectively. For the definition of these topologies, see Skorohod \cite{Skorokhod1956} or
\cite{Bill99}, \cite{LouhichiRio2011},
\cite{ResnickVdB2000}.  Denote $|\mu|_{2+\delta}(u) = \E |\gamma(u)|^{2+\delta} $
the absolute $(2+\delta)$-moment of the seed process.

\subsection{Gaussian scenario (Case 1)}\label{31}


\begin{theo} \label{thmgauss}
\
\begin{itemize}
\item[(i)]  Assume $\mu(u) = \E \gamma(u) = 0$,   \eqref{g0},  \eqref{G11}, \eqref{g4} and
\begin{equation}\label{G2}
|\mu|_{2+\delta}(u) 
= {\cal O}(|u|^{\frac{2+\delta}2}), \qquad (u \to 0, \, \exists \, \delta >0).
\end{equation}
Then
\begin{equation}\label{cltD}
\frac 1{n^H} S_{[nt]} \ \to_{{\mathcal D}(J_1)} \ \sqrt{ c_2} \,B_H(t),  \qquad H = \frac{3-\alpha}2
\end{equation}
where $B_H$ is fractional Brownian motion with variance $\E B^2_H (t) = t^{2H} $ and
$c_2$ is defined in  \eqref{g5}.

\item[(ii)] Assume $\mu(u) = \E \gamma(u) = 0$,  \eqref{rhobdd}, \eqref{G4},  \eqref{G2}
and $\sigma^2  =
\sum_{k\in \Z} \cov   (X_0, X_k) \ne 0$. Then
\begin{equation}\label{cltD2}
\frac 1{\sqrt{n}} S_{[nt]} \ \to_{f.d.d.} \  \sigma \,B(t),
\end{equation}
where $B$ is a Brownian motion with variance $\E B^2(t) = t$.\\
In addition, if \ $\sum_{k=1}^\infty \sqrt{|a_k|} < \infty $, then the finite dimensional convergence
in \eqref{cltD2} can be replaced by 
$\to_{{\mathcal D}(J_1)}$.

\item[(iii)] Assume the same conditions as in (ii) except that \eqref{G2} is replaced by
\begin{equation}\label{mu0}
|\mu|_{2+\delta}(u) =  {\cal O}(u) \quad (u \to 0) \quad \text{and} \quad \sum_{j=0}^\infty |a_j|^{\frac1{2+\delta}} < \infty
\end{equation}
for some $\delta >0$. Then all statements in part (ii) remain valid.

\end{itemize}
\end{theo}


\noi {\it Proof.} (i)
Consider the convergence of  one-dimensional distributions:
\begin{equation}\label{clt}
\frac 1{\sqrt{n^{3-\alpha}}} S_{n} \ \to_{law} \   {\cal N}(0,c_2).
\end{equation}
In view of \eqref{g5} and Lemma \ref{decomplem},  relation \eqref{clt} follows by
Lindeberg's theorem provided
\begin{equation} \label{lind}
L_n    \ :=  \ \sum_{s=-\infty}^n   \E |Z_{s,n}|^{2+\delta} \  = \  o\big(n^{\frac{(3-\alpha)(2+\delta)}2}\big).
\end{equation}
By Minkowski's inequality and assumptions \eqref{G12}  and \eqref{G2} 
we obtain
\begin{eqnarray}\label{Zsn}
\E |Z_{s,n}|^{2+\delta}&\le&\left( \sum_{k=1\vee s}^n  (\E |\gamma(a_{k-s})|^{2+\delta})^{\frac1{2+\delta}} \right)^{2+\delta} \\
&\le&C\left( \sum_{k=1\vee s}^n  |a_{k-s}|^{\frac12} \right)^{2+\delta} \
\le\  C\left( \sum_{k=1\vee s}^n  |k-s|_+^{-\frac\alpha2} \right)^{2+\delta} \nn
\end{eqnarray}
and therefore $L_n  \le  C(L^-_n + L^+_n)$, where
\begin{eqnarray*}
L^-_n&=&\sum_{s=-\infty}^0  \left( \sum_{k=1}^n  |k-s|_+^{-\frac\alpha2} \right)^{2+\delta}
\ = \  \sum_{s=0}^\infty  \left( \sum_{k=1}^n  (k+s)^{-\frac\alpha2} \right)^{2+\delta}, \\
L^+_n&=&\sum_{s=1}^n \left( \sum_{k=1}^n  k^{-\frac\alpha2} \right)^{2+\delta} \ = \ n \left( \sum_{k=1}^n  k^{-\frac\alpha2} \right)^{2+\delta}.
\end{eqnarray*}
Here, $L^+_n  = {\cal O}\left( n \big( n^{1-\frac \alpha2})^{2+\delta}\right)  =   o\big(n^{\frac{(3-\alpha)(2+\delta)}2}\big)$.
The same relation for $L^-_n $ follows from
\begin{eqnarray*}
L^-_n&\le&\int_{0}^\infty \d x \left( \int_0^n  (x+y)^{-\frac\alpha2} \d x \right)^{2+\delta}\ = \ c n \left( n^{1- \frac\alpha2}\right)^{2+\delta}, \quad \mbox{ with}\\
c& =&  \int_{0}^\infty \d x \left( \int_0^1 (x+y)^{-\frac\alpha2} \d x \right)^{2+\delta} < \infty.
\end{eqnarray*}
This proves \eqref{lind} and the one-dimensional convergence in \eqref{clt}.  Finite-dimensional convergence
in \eqref{cltD} follows similarly using Cram\'er-Wold device. Finally, the tightness in ${\mathcal D}(J_1)$ of the partial sums process
in \eqref{cltD} follows by Kolmogorov's criterion and from property  \eqref{g5} (see, e.g.
\cite{gir2012}, proposition~4.2.2). This proves part (i).

\smallskip

\noi (ii) Again, it suffices to prove  the convergence of  one-dimensional distributions:
\begin{equation}\label{clt2}
n^{-1/2} S_{n} \ \to_{law} \   {\cal N}(0,\sigma^2).
\end{equation}
By writing $S_n$ as in \eqref{Zn} and using Lindeberg's theorem relation \eqref{clt2} follows from
\begin{equation} \label{lind2}
L_n    \ =  \ \sum_{s=-\infty}^n   \E |Z_{s,n}|^{2+\delta} \  = \  o\big(n^{\frac{2+\delta}2}\big).
\end{equation}
Using Minkowski's inequality and assumptions \eqref{G2} and \eqref{G4} similarly as in part (i) we obtain
\begin{eqnarray} \label{Zbdd}
\E |Z_{s,n}|^{2+\delta}&\le&C\Big( \sum_{k=1\vee s}^n  |a_{k-s}|^{\frac12} \Big)^{2+\delta} \\
&\le&C\Big( \sum_{k=1\vee s}^n  |(k-s)a_{k-s}| \Big)^{\frac{2+\delta}2} \Big( \sum_{k=1\vee s}^n  (k-s)^{-1} \Big)^{\frac{2+\delta}2} \nonumber \\
&\le&C\Big( \sum_{k=1\vee s}^n  (k-s)^{-1} \Big)^{\frac{2+\delta}2}.
\end{eqnarray}
and hence 
\begin{eqnarray*}
\sum_{s=-n}^n   \E |Z_{s,n}|^{2+\delta}&\le&C n  (\log n)^{\frac{2+\delta}2} \ = \  o \big(n^{\frac{2+\delta}2}\big), \\
\sum_{s=-\infty}^{-n}   \E |Z_{s,n}|^{2+\delta}
&\le&C \sum_{s=n}^\infty \Big( \sum_{k=1}^n \frac1{k+s} \Big)^{\frac{2+\delta}2} \ \le \
C \sum_{s=n}^\infty ( n s^{-1} )^{\frac{2+\delta}2}  \ \le \ C n  \ =  \  o \big(n^{\frac{2+\delta}2}\big),
\end{eqnarray*}
proving \eqref{lind2} and \eqref{clt2}.
To show the last statement of (ii) (the
tightness in $D[0,1]$), it suffices to prove the bound
\begin{equation}\label{Sp}
\E |S_n|^{2+\delta} \ \le \ C n^{\frac {2+\delta}{2}},
\end{equation}
see (\cite{gir2012}, proposition~4.4.4). By Rosenthal's inequality,
$$
\E |S_n|^{2+\delta} \ \le \ C \Big( \sum_{s=-\infty}^n   (\E |Z_{s,n}|^{2+\delta})^{\frac{2}{2+\delta}}  \Big)^{\frac{2+\delta}{2}}.
$$
Using \eqref{Zbdd} and $\sum_{k=1}^\infty |a_k|^{\frac12} < \infty $, we get
$\max_{|s| \le n} \E |Z_{s,n}|^{2+\delta} < C $ and
\begin{eqnarray}
 \sum_{s=-\infty}^{-n}   (\E |Z_{s,n}|^{2+\delta})^{\frac2{2+\delta}}
&\le&C\sum_{s=n}^\infty \Big(\sum_{k=1}^n |a_{k+s}|^{\frac12} \Big)^2 \nn \\
&\le&C\sum_{k_1, k_2=1}^n \sum_{s=n}^\infty |a_{k_1+s}|^{\frac12} |a_{k_2+s}|^{\frac 12} \  \le \ Cn. \label{2nd}
\end{eqnarray}
This proves \eqref{Sp} and part (ii), too.

\smallskip


\noi (iii) Similarly as in \eqref{Zsn} and using \eqref{mu0} we get  
$$ 
\E |Z_{s,n}|^{2+\delta} \le C\Big( \sum_{k=1\vee s}^n  |a_{k-s}|^{\frac{1}{2+\delta}} \Big)^{2+\delta} \le C\sum_{k=1\vee s}^n  |a_{k-s}|^{\frac{1}{2+\delta}} \le C $$
for any $-\infty < s \le n$  and hence
\begin{eqnarray*}\label{Zsn0}
\sum_{s=-\infty}^{-n} \E |Z_{s,n}|^{2+\delta}&\le&C\sum_{s=n}^\infty  \sum_{k=1}^n  |a_{k+s}|^{\frac{1}{2+\delta}}  \ \le \ Cn, \\
\sum_{s=-\infty}^{-n} (\E |Z_{s,n}|^{2+\delta})^{\frac{2}{2+\delta}}&\le&C
\sum_{s=n}^\infty \big( \sum_{k=1}^n  |a_{k+s}|^{\frac{1}{2+\delta}} \big)^2 \ \le \ Cn
\end{eqnarray*}
as in \eqref{2nd}. Hence, \eqref{lind2} and \eqref{Sp} follow, proving part (iii) and completing the proof of
Theorem \ref{thmgauss}.
\hfill $\Box$

\begin{Rem} 
{\rm
The crucial condition for Gaussian partial sums limit  under long-range
dependence assumption \eqref{g4}  in Theorem \ref{thmgauss} (i)
is
\eqref{G2}. Clearly this condition is satisfied for the Brownian motion
$\gamma (u) = B(u) $, in which case $|\mu|_{2+\delta}(u) = \E |B(u)|^{2+\delta} = |u|^{\frac{2+\delta}2} \E |B(1)|^{2+\delta}  $.
On the other hand, condition \eqref{G2} is not satisfied for most jump processes.
Particularly, if $\gamma (u) = P(u) - u,  u \ge 0$ is a centered Poisson process with intensity
$\E P(u) = u$, then 
$$
|\mu|_{2+\delta}(u) = u \e^{-u} |1-u|^{2+\delta} + {\cal O}(u^{2+ \delta} + u^2)
\sim u  \quad (u \to 0)
$$
and \eqref{G2} fails, but the first condition in \eqref{mu0} is satisfied. In particular,
in the case of Poisson seed process, the trawl process satisfies Donsker's theorem if
the trawl decays fast enough so that  \eqref{mu0} holds.
}
\end{Rem}

\medskip

\noi
Let us present  further examples of seed processes satisfying the conditions in Theorem \ref{thmgauss}.

\smallskip

\begin{Ex}[Geometric centered Brownian motion] \label{Ex4} {\rm
Set  $\gamma(u) = \e^{ B(u) - u/2} -1, u \ge 0$, where
$B$ is a standard Brownian motion as above. We have   $\E \gamma(u) = 0$ and (if $u\le v$)
\begin{eqnarray*}
\rho(u,v)&=&\E\exp\{ B(u) + B(v) - \frac{u+v}2\} -1 \\
&=&\exp\Big\{ (\frac12 \E (B(u) + B(v))^2 - \frac{u+v}2\Big\} - 1 \\
&=&\exp\Big\{ (\frac12 (u + v + 2 u)   -  \frac{u+v}2\Big\} - 1   \\
&=&\e^u -1  = u\wedge v + O \big( (u\wedge v)^2\big), \quad  u\wedge v \to 0.
\end{eqnarray*}
Therefore conditions \eqref{g0}, \eqref{G11} are satisfied.
We also have by Taylor's expansion that
$|\mu|_{4}(u) = \E \left| \e^{ B(u) - u/2} -1\right|^{4}
= \e^{6u} - 4 \e^{3u} + 6 \e^u - 3= {\cal O}(u^2), \ u \to 0 $
so that \eqref{G2} is satisfied with $\delta = 2 $.
}
\end{Ex}

\begin{Ex} [Diffusion process]\label{Ex5}
{\rm
$$
\gamma(u) = \int_0^u b(v) \d B(v)
$$
with $B $ a Brownian motion, and $(b(v))_{v \ge 0}$  a random predictable process with
$\lim_{v\to 0}\E b^2(v)=  C >0 $. Then $g(u) = \int_0^u \E b^2(v) \d v \sim  C u \, (u \to 0) $ and
$\rho(u,v)= g(u), \ 0 \le u \le v $
so that conditions  \eqref{g0} and \eqref{G11} are satisfied. Moreover, if $\E |b(v)|^{2+\delta} \le C$ then
by the moment inequality for Brownian integrals (see, e.g. \cite{kwa1992}, theorem~9.9.2)
\begin{eqnarray*}
|\mu|_{2+\delta}(u)&\le&C\E \left(\int_0^u   b^2(v) \d v\right)^{\frac{2+\delta}2} \\
&\le&C \left(\int_0^u  \E |b(v)|^{2+\delta} \d v \right) \left( \int_0^u 1 \, \d v \right)^{\frac{2+ \delta}2 - 1} \
\le\  Cu^{\frac{2+\delta}2},
\end{eqnarray*}
hence assumption \eqref{G2} holds, too.
 }
\end{Ex}

\subsection{Stable scenario (Case 2)}\label{32}

We assume now that seed process $\gamma = \{\gamma(u), u \ge 0 \}$ is a piece\-wise constant nondecreasing process starting at $\gamma(0)=0$
with unit jumps at points $0 = \tau_0 < \tau_1 < \tau_2 < \cdots $:
\begin{equation} \label{Z0}
\gamma(u)= \sum_{k=0}^\infty k\cdot \1(\tau_k \le u < \tau_{k+1})  
\end{equation}
and such 
that the distribution of the first jump-point $\tau_1 >0$ has a bounded probability density
$\theta(u)$:
\begin{equation}\label{tau1}
\P(0 < \tau_1 \le u) = \int_0^u \theta(y) \d y  \qquad \text{and} \qquad \lim_{u \to 0} \theta(u) = 1.
\end{equation}
Moreover, we shall assume that there exists $\delta >2(\alpha-1)$ such that
\begin{eqnarray}
&&\E \gamma(u)^{2+ \delta}\ <\  \infty, \qquad \qquad \qquad \forall \  u >0,  \label{Z1} \\
&&\E \gamma(u)^2 \1 (\tau_2 \le u)\ =\  {\cal O}(u^2), \qquad u \to 0.    \label{Z2}
\end{eqnarray}
\begin{Rem}
{\rm The second condition in \eqref{tau1} can be replaced by $\lim_{u \to 0} \theta(u) = C >0$ without loss of generality.
Conditions \eqref{tau1}-\eqref{Z3}  are very general and there are satisfied by many jump processes
$\gamma$ as this was sketched in the introduction.  As shown below, these conditions also imply the conditions on $\gamma $ in Proposition \ref{propLM}.
\\
Remark that $(\tau_1\le u)=(\gamma(u)\ge 1)$, $(\tau_2\le u)=(\gamma(u)\ge2)$ and therefore an alternative way to set  condition \eqref{Z2} is $\E\gamma^2(u){\1}(\gamma(u)>1)= {\cal O}(u^2)$, as $ u \to 0$.
}
\end{Rem}

\begin{prop} For the seed process $\gamma$ in \eqref{Z0},
conditions  \eqref{tau1}-\eqref{Z2} imply the
assumptions \eqref{mug} and \eqref{g0} of Proposition \ref{propLM}~(i). In addition, if
\begin{eqnarray}
&&\E \gamma(v) \1 (\tau_1 \le u, \tau_2 \le v)\ =\  o(u), \qquad 0 \le u \le v  \to 0,    \label{Z3}
\end{eqnarray}
then  \eqref{G11} is satisfied.
\end{prop}

\noi {\it Proof.} From \eqref{Z0}  we have
\begin{equation}\label{zetabdd}
\1 (\tau_1 \le u) \ \le \ \gamma(u) \ \le \ \1 (\tau_1 \le u) + \gamma(u) \1 (\tau_2 \le u)
\end{equation}
and hence
\begin{equation*}
\P (\tau_1 \le u) \ \le \ \mu(u) \ \le \ \P(\1 (\tau_1 \le u) + \E \gamma(u) \1 (\tau_2 \le u)
\end{equation*}
From  \eqref{tau1},
$\P(0 < \tau_1 \le u) = u (1 + o(1)) $ and from \eqref{Z2},
$$
\E \gamma(u) \1 (\tau_2 \le u) \le \E \gamma^2(u) \1 (\tau_2 \le u)   = {\cal O}(u^2).
$$
Therefore,
\begin{equation} \label{mu1}
\mu (u)\  = \ u (1 + o(1)) + {\cal O}(u^2) \ = \   u (1 + o(1)) \qquad (u \to 0).
\end{equation}
Similarly, for the second moment $\mu_2(u) = \E \gamma^2(u)$ from \eqref{zetabdd}, \eqref{tau1}, \eqref{Z2} we obtain
$$
\P (\tau_1 \le u) \ \le \ \mu_2(u) \ \le \ \P(\1 (\tau_1 \le u) + 2\E \gamma(u) \1 (\tau_2 \le u) +
\E \gamma^2(u) \1(\tau_2 \le u),
$$
implying
$\mu_2(u)  =  u (1 + o(1)) + {\cal O}(u^2)  =  u (1 + o(1)) \ (u \to 0)$
and
\begin{equation}\label{mu3}
g(u)\  = \ \mu_2(u) - \mu^2(u)
 \ = \  u (1 + o(1)) \qquad (u \to 0).
\end{equation}
Clearly, \eqref{mu1}  and \eqref{mu3} imply \eqref{mug} and \eqref{g0}.
Consider assumption \eqref{G11}. Since
$$
\rho (u,v) = \E \gamma(u) \gamma(v) - \mu(u) \mu (v)
=  \E \gamma(u) \gamma(v) - u v (1 + o(1)) = \E \gamma(u) \gamma(v) + o(u \wedge v),
$$
as $ 0 < u \le v \to 0$, condition
\eqref{G11} follows from
\begin{equation} \label{Ebdd}
\E \gamma (u)  \gamma(v)  = u(1 + o(1)), \qquad   0 < u \le v \to 0.
\end{equation}
From  \eqref{zetabdd} for $0< u \le v $ we obtain
\begin{eqnarray*}
\P (\tau_1 \le u)&\le&\E \gamma (u)  \gamma(v)  \\
&\le&\P(\tau_1 \le u)
+\E \gamma(u) \1 (\tau_2 \le u) +
\E \gamma(v) \1(\tau_1 \le u, \tau_2 \le v) + \E \gamma(u) \gamma(v) \1 (\tau_2 \le u)
\end{eqnarray*}
where
$\E \gamma(u) \gamma(v) \1 (\tau_2 \le u) \le (\E \gamma^2(u) \1 (\tau_2\le u))^{\frac12}
(\E \gamma^2 (v))^{\frac12}  \le  C u (\E \gamma^2 (v))^{\frac12} $
and $\E \gamma^2 (v) = \mu_2(v) =  {\cal O}(v),  $
see \eqref{tau1}, \eqref{Z2}.
Hence and from \eqref{Z3} we have that
$$
\E \gamma(u) \1 (\tau_2 \le u) +
\E \gamma(v) \1(\tau_1 \le u, \tau_2 \le v) + \E \gamma(u) \gamma(v) \1 (\tau_2 \le u) = o(u)
$$
implying \eqref{Ebdd} and  \eqref{G11}, too. \hfill $\Box$

\begin{theo} \label{thmstable} Assume that $a_j \ge 0$  satisfy the regular decay condition in \eqref{g4} with exponent
$1< \alpha < 2 $ and that the seed process in \eqref{Z0} satisfies conditions \eqref{tau1}-\eqref{Z2}.
Then
\begin{equation}\label{cltL}
n^{-\frac1\alpha} (S_{[nt]} - \E S_{[nt]}) \ \to_{f.d.d.} \   L_\alpha(t),
\end{equation}
where $L_\alpha (t), t \ge 0$ is a homogeneous $\alpha$-stable L\'evy process with
characteristic function
\begin{equation}\label{chfL}
\E \e^{ \i z L_\alpha (t)} \ = \  \exp \left\{ - t |z|^\alpha \frac{c_0 \Gamma (2-\alpha)}{1-\alpha}
\left( \cos (\pi \frac\alpha2) - \i\cdot {\rm sgn}(z) \sin (\pi \frac\alpha 2) \right) \right\}, \quad z \in \R.
\end{equation}

\end{theo}

\noi {\it Proof.}
Denote
\begin{equation}\label{Zrv}
Z \ = \ \sum_{j=0}^\infty \gamma (a_j), \quad Z^* \ = \  \sum_{j= 0}^\infty \1(\gamma(a_j) \ge 1) \ = \  \# \{  j\ge 0: a_j \ge \tau_1 \},
\quad  Z^{**} \  = \ Z - Z^*.
\end{equation}
Then $Z \ge Z^* \ge 0$ and the series for $Z$ in \eqref{Zrv} converges a.s. in view of \eqref{mu1}
and has finite mean:
$$\E Z = \sum_{j=0}^\infty \mu (a_j) \le C \sum_{j=0}^\infty a_j < \infty.$$
We shall prove that the tail d.f. of r.v. $Z$ decays regularly with exponent $\alpha \in (1,2)$:
\begin{equation} \label{Ztail}
\P (Z > y) \ =  \  c_0 \, y^{-\alpha}(1 + o(1)), \qquad \text{as} \ \  y \to \infty.
\end{equation}
Relation \eqref{Ztail}  follows from \eqref{Zrv} and
\begin{equation} \label{Z1tail}
\P (Z^* > y) \ =  \  c_0 \, y^{-\alpha}(1 + o(1))   \quad \text{and}  \quad
\P (Z^{**} > y) \ =  \  o(y^{-\alpha}), \quad \text{as} \ \  y \to \infty,
\end{equation}
\noi
Consider the first relation in  \eqref{Z1tail}. Since $\P (Z^* > k-1) \ge \P (Z^* > y) \ge  \P (Z^* > k) $ when $k-1 \le y \le k $, it suffices
to show \eqref{Z1tail} for $y = k-1 $, or the probability $\P (Z^* \ge k), k \in \N_+$.
As noted in the proof of Proposition \ref{propLM},
for any $\epsilon >0$  there exists $j_0 > 0$ such that
$c_0 (1-\epsilon) j^{-\alpha} < a_j < c_0 (1+ \epsilon) j^{-\alpha}, \, \forall \, j\ge j_0$. Clearly, for any
$k \ge 1 $ we have
$ \P (Z_- \ge k+j_0) \le \P (Z^* \ge k) \le \P(Z_+ \ge k-j_0)$, where
\begin{eqnarray*}
&Z_+ = \sum_{j=j_0}^\infty \1 (\tau_1 \le c_0(1+\epsilon)j^{-\alpha}) = \# \{j \ge j_0:  \tau_1 \le c_0(1+\epsilon)j^{-\alpha}\}, \\
&Z_- = \sum_{j=j_0}^\infty \1 (\tau_1 \le c_0(1-\epsilon)j^{-\alpha}) =  \# \{j \ge j_0:  \tau_1 \le c_0(1+\epsilon)j^{-\alpha}\}.
\end{eqnarray*}
According to 
\eqref{tau1}, as $k \to \infty$,
\begin{equation*}
\P (Z_+ \ge k-j_0) \ = \  \P (\tau_1 < c_0 (1+\epsilon) k^{-\alpha} )\ = \  \int_0^{c_0 (1+\epsilon) k^{-\alpha}} \theta (y) \d y
\  \sim \ c_0(1+\epsilon) \, k^{-\alpha}
\end{equation*}
and, similarly,
\begin{equation*}
\P (Z_- \ge k+j_0) \ = \  \P (\tau_1 < c_0 (1-\epsilon) (k + 2j_0 -1)^{-\alpha} )
 \  \sim \ c_0(1-\epsilon) \, k^{-\alpha}.
\end{equation*}
Therefore, $c_0 (1- \epsilon)  \le    \liminf k^{\alpha} \P (Z^* \ge k) \le  \limsup k^{\alpha} \P (Z^* \ge k) \le c_0 (1+ \epsilon)$,
where $\epsilon >0$ is arbitrary small,
proving the first relation in  \eqref{Z1tail}. To prove  the second relation in  \eqref{Z1tail}, note
$Z^{**} \le \sum_{j=0}^\infty  \gamma(a_j) \1(a_j \ge \tau_2) $ and then 
by  \eqref{Z2} and
Minkowski's inequality we obtain
\begin{eqnarray*}
\E^{\frac12} (Z^{**})^2
&\le&\sum_{j=0}^\infty \big( \E \gamma^2(a_j) \1(a_j \ge \tau_2)\big)^{\frac12} \
\le\  C \sum_{j=0}^\infty |a_j|  \  < \ \infty
\label{Z5}
\end{eqnarray*}
proving  the second relation in  \eqref{Z1tail} and hence \eqref{Ztail} as well. In turn,
\eqref{Ztail} implies that the distribution of r.v. $Z$ belongs to the domain
of attraction of asymmetric $\alpha$-stable  law, viz.,
\begin{equation}\label{cltZ}
n^{-\frac 1\alpha}  \sum_{k=1}^{[nt]} (Z_k - \E Z_k) \ \to_{f.d.d.} \   L_\alpha(t),
\end{equation}
where
$Z_k = \sum_{j=0}^\infty \gamma_k(a_{j}), \, k \in \Z$
are i.i.d. copies of r.v. $Z$ in \eqref{Zrv} and
$L_\alpha $ is the $\alpha$-stable L\'evy process in \eqref{cltD}-\eqref{chfL}.
See e.g. (\cite{ibr1971}, theorem~2.6.7).
\\
Relation \eqref{cltL} follows from \eqref{cltZ} if we show
that the partial sums process in \eqref{cltL} can be approximated by the partial sums process
in \eqref{cltZ}, in the sense that
\begin{equation}\label{ES}
\E |S_n - \widetilde S_n| \ = \  o(n^{\frac1\alpha}), \qquad \text{where} \quad \widetilde S_n
= \sum_{k=1}^n Z_k.
\end{equation}
We have $\widetilde S_n - S_n  = R'_n - R''_n$, where
$$
R'_n = \sum_{1 \le s \le n} \sum_{j > n-s} \gamma_s (a_j), \qquad
R''_n = \sum_{s \le 0} \sum_{1\le k \le n} \gamma_s (a_{k-s}), $$
then $ R'_n \ge 0, R''_n \ge 0$.
\\
Using  \eqref{mu1} and \eqref{g4} we obtain
\begin{eqnarray*}
\E R'_n&=&\sum_{1 \le s \le n} \sum_{j > n-s} \E \gamma_s (a_j) \ = \ \sum_{1 \le s \le n} \sum_{j > n-s} \mu(a_j) \\
&\le&C \sum_{1 \le s \le n} \sum_{j > n-s} j^{-\alpha} \ = \ {\cal O}(n^{2-\alpha}),    \\
\E R''_n &=&\sum_{s \le 0} \sum_{1 \le k \le n} \E \gamma_s (a_{k-s}) \ = \  \sum_{s \le 0} \sum_{1 \le k\le n} \mu (a_{k-s}) \\
&=&\sum_{s \ge 0} \sum_{1\le k \le n} (k+s)^{-\alpha}   \ = \ {\cal O}(n^{2-\alpha}),
\end{eqnarray*}
implying \eqref{ES} since $ 2- \alpha < 1/\alpha$ for $1 < \alpha < 2 $. Theorem
\ref{thmstable} is proved.  \hfill  $\Box$

\medskip

\begin{Ex} [Jump processes and the assumptions of Theorem  \ref{thmstable}]
\label{Rem1}{\rm  
For such a jump process  $(\gamma(u)=k)=(\tau_k\le u<\tau_{k+1})$.
Conditions in \eqref{tau1}-\eqref{Z2} on the seed process
$\{\gamma(u), u \ge 0\} $ in Theorem \ref{thmstable} are rather weak and essentially involve the
distribution of the first jump-point $\tau_1$ provided the second jump $\tau_2$  cannot occur
very fast after $\tau_1$. Particularly, 
\begin{itemize}
\item
 The Bernoulli process  is very simple: in this  case $\tau_2 = \infty $ thus $\P(\tau_2\le u)= \P(\tau_2\le v)
= 0$,

\item The Poisson  process in Example \ref{ex3}. Indeed, for $\gamma (u) = P(u) $
\eqref{Z2} holds  since
\begin{eqnarray*}
&\E \gamma(u)^2 \1 (\tau_2 \le u)\
\ =\ \E \gamma(u)^2 - \P(\gamma(u) =1)\  = \ u + u^2  - u \e^{-u} \ = \  {\cal O}(u^2).
\end{eqnarray*}
Verification of \eqref{Z3} for  $\gamma (u) = P(u) $ is slightly more involved, as follows.  Let $b(u,v) = \E \gamma(v) \1 (\tau_1 \le u, \tau_2 \le v)
= b_1(u,v) + b_2(u,v)$, where $b_1(u,v) = \E \gamma(v) \1 ( \tau_2 \le u) \le \P^{2/3}(\tau_2 \le u) \E^{1/3} \gamma^3(v) =
{\cal O}(u^{4/3}) = o(u) $ since  $\P(\tau_2 \le u) = \P(\gamma(u) \ge 2) = {\cal O}(u^2), u \to 0$. Next, since $\tau_2 >u$ implies $\gamma(u) = 1$ so
$b_2(u,v) = \E \gamma(v) \1 (\tau_1 \le u, u< \tau_2 \le v) =  \P (\tau_1 \le u, u< \tau_2 \le v)
+  \E  (\gamma(v) - \gamma(u)) \1 (\gamma(u)=1, \gamma(v) \ge 1), $ where
$ \P (\tau_1 \le u, u< \tau_2 \le v) = \P(\gamma(u) =1, \gamma(v) -\gamma(u) \ge 1) = {\cal O}\big(u(v-u)\big) = o(u)$ and,
similarly  $\E  (\gamma(v) - \gamma(u)) \1 (\gamma(u)=1, \gamma(v) \ge 1)
= \P (\gamma(u) = 1) \E \gamma (v-u) =  {\cal O}(u(v-u)) = o(u),  \, 0< u \le v \to 0$, proving
\eqref{Z3}. 
\item
Other examples of jump processes satisfying \eqref{tau1}-\eqref{Z2} include
mixed Poisson processes (Example \ref{ex3}) and
renewal process with independent intervals $\tau_1$ and $\tau_2 - \tau_1 $
and $\P (\tau_2 - \tau_1 \le x) = {\cal O}(x) $
since
$$
\P(\tau_2 \le  u) = \int_0^u \theta(y) \P (\tau_2 - \tau_1\le u-y)   \d y
\le C \int_0^u  (u-y)   \d y  = {\cal O}(u^2)
$$
as in the Poisson case.
The same  conditions also holds for mixed Poisson processes driven by some random variable $\zeta>0$ (Example \ref{ex3}).
(thus again the case of negative binomials fits our result as sketched in \cite{fok}).
\end{itemize}

}
\end{Ex}


\noi
We note that the functional convergence in  \eqref{cltL} is open and may
not hold in the  $J_1$-topology.
At the cost of additional structure we can prove the convergence in Skorohod's $M_1$-topology.
For definitions and properties related to association of random variables we refer to  \cite{EPW1967}.

\begin{theo} \label{thmfunct} Suppose that all assumptions of Theorem \ref{thmstable} hold. If the jump random variables $\tau_1, \tau_2, \ldots $ are {\em associated} (in particular, if they are sums of independent positive random variables) then the finite-dimensional convergence (\ref{cltL})
can be strengthened to
\begin{equation}\label{cltLF}
n^{-\frac1\alpha} (S_{[nt]} - \E S_{[nt]}) \ \to_{{\mathcal D}(M_1)} \   L_\alpha(t),
\end{equation}

\end{theo}

\noi {\it Proof.}  By (\cite{LouhichiRio2011}, theorem 1) it suffices to verify that $X_1, X_2, X_3, \ldots $ are associated random variables. By (\cite{EPW1967}, property $P_5$) it is enough to check association of
\[
\sum_{j=0}^N \gamma_{1-j}(a_j), \quad\sum_{j=0}^N \gamma_{2-j}(a_j), \ldots,\quad \sum_{j=0}^N \gamma_{k-j}(a_j),
\]
for each $N\in \N$ and $k\in \N$, where $\gamma_j(\cdot)$ are independent copies of (\ref{Z0}). This in turn is implied by (\cite{EPW1967}, properties $P_4$ and $P_2$), provided the family
$
\gamma(a_1), \gamma(a_2), \gamma(a_3),\ldots,
$
is associated. But
\[ \gamma(u)= \sum_{k=0}^\infty k\cdot \1(\tau_k \le u < \tau_{k+1}) = \sum_{k=1}^{\infty} \1(\tau_k \le u),\]
and by arguments already presented above it is enough to prove association of random variables
\begin{equation}
\label{assoarray}
\begin{array}{cccc}
\1(\tau_1 \le a_0), & \1(\tau_2 \le a_0), & \cdots & \1(\tau_k \le a_0), \\
\1(\tau_1 \le a_1), & \1(\tau_2 \le a_1), & \cdots & \1(\tau_k \le a_1), \\
\vdots & \vdots & \ddots & \vdots \\
\1(\tau_1 \le a_N), &\1(\tau_2 \le a_N), & \cdots & \1(\tau_k \le a_N),
\end{array}
\end{equation}
for each $N\in \N$ and $k\in \N$. Let us notice that
$\1(\tau_j \le a_m),  = 1 -  \1(\tau_j > a_m)$ and that functions $f_m(x) = \1(x > a_m)$ are nondecreasing. Therefore, if $\tau_1, \tau_2, \ldots $ are associated, then also the family $\{ f_m(\tau_j) \,;\, j,m \in \N\}$ is associated. By (\cite{EPW1967}, property $BP_1$) array (\ref{assoarray}) is associated as well.
\hfill $\Box$

\begin{Rem}\label{plusminus} {\rm
Let us consider two families $\{\gamma_j^+\}$ and $\{\gamma_j^-\}$ of processes of the form (\ref{Z0}). Consider stationary processes $X_1^+, X_2^+, \ldots $, and $X_1^-, X_2^-, \ldots $, each built according to the recipe (\ref{Xk}), and the corresponding partial sum processes $S_{[nt]}^+$ and $S_{[nt]}^-$.

If $\{\gamma_j^+\}$ and $\{\gamma_j^-\}$ are independent and
\begin{equation}\label{fddsing}
 n^{-\frac1\alpha} (S^+_{[nt]} - \E S^+_{[nt]}) \ \to_{f.d.d.} \   L^+_\alpha(t), \quad  n^{-\frac1\alpha} (S^-_{[nt]} - \E S^-_{[nt]}) \ \to_{f.d.d.} \   L^-_\alpha(t),
\end{equation}
then also
\begin{equation}\label{fddsum}
 n^{-\frac1\alpha} \Big(\big(S^+_{[nt]} - S^-_{[nt]}\big) - \E \big(S^+_{[nt]} - S^-_{[nt]}\big)\Big) \ \to_{f.d.d.} \   L_\alpha(t),
\end{equation}
where $L_{\alpha} \sim  L^+  - L^-$ with independent $L^+ \sim L^+_\alpha$ and $ L^- \sim L^-_\alpha$.
\medskip

\noindent
In particular, if $\gamma^+_j$ and $\gamma^-_j$ are identically distributed, then the resulting trawl process is centered and the limiting L\'evy process is
{\em symmetric}. This is the case if e.g. $\gamma^\pm$ are both homogeneous Poisson processes with identical intensities or Bernoulli processes $\gamma^\pm(u)={\1}(U^\pm\le u)$ for independent uniform rvs, $U^\pm$.
}
\end{Rem}

\begin{Rem}
{\rm
As the example of an ordinary moving average with summable coefficients shows, (\ref{fddsing}) may imply (\ref{fddsum}) without the assumption
of independence of $S^+_{[nt]}$ and $S^-_{[nt]}$ (see e.g. (\cite{BJL2016}, corollary 2.2)).
In the functional limit theorem given below we follow this general approach and obtain the functional convergence in the non-Skorohodian $S$ topology (see \cite{Jakubowski1997}).
We shall denote by $\to_{{\mathcal D}(S)}$
the weak convergence  in the Skorohod space $D[0,1]$ equipped with the $S$ topology).
}
\end{Rem}

\begin{cor} \label{corStop} In the framework of Remark \ref{plusminus}, suppose  that both $S^+_{[nt]}$ and $S^-_{[nt]}$ satisfy all assumptions of
Theorem \ref{thmfunct}, so that
\begin{equation}\label{functem}
 n^{-\frac1\alpha} (S^+_{[nt]} - \E S^+_{[nt]}) \ \to_{\mathcal{D}(M_1)} \   L^+_\alpha(t), \quad  n^{-\frac1\alpha} (S^-_{[nt]} - \E S^-_{[nt]}) \ \to_{\mathcal{D}(M_1)} \   L^-_\alpha(t),
\end{equation}
for some $\alpha$-stable L\'evy motions  $L^+_\alpha$ and  $L^-_\alpha$.

If for some c\`adl\`ag stochastic process $K$ we have
\begin{equation}\label{fddsumS}
 n^{-\frac1\alpha} \Big(\big(S^+_{[nt]} - S^-_{[nt]}\big) - \E \big(S^+_{[nt]} - S^-_{[nt]}\big)\Big) \ \to_{f.d.d.} \   K(t),
\end{equation}
then
\[ n^{-\frac1\alpha} \Big(\big(S^+_{[nt]} - S^-_{[nt]}\big) - \E \big(S^+_{[nt]} - S^-_{[nt]}\big)\Big) \ \to_{{\mathcal D}(S)}
\   K(t).
\]

\end{cor}

\noi {\it Proof.} By (\cite{BJL2016}, theorem 3.13) (\ref{functem}) implies
the uniform $S$-tightness of the corresponding processes. The proof of  (\cite{BJL2016}, proposition 3.16) gives the uniform $S$-tightness of the differences. A direct application of  (\cite{BJL2016}, proposition 3.3) concludes the proof.
\hfill $\Box$.

\paragraph{Acknowledgements.}
This study begun with a question from Wilfredo Palma (Santiago de Chile)  to the first author: {\it how to define LRD integer valued models?} We wish to thank him for considering this problem.
\\
This work has been developed within the MME-DII center of excellence (ANR-11-LABEX-0023-01) and was partially supported by CNPq-Brazil.
\\
We also thank the Universities UFRGS  (Porto Alegre) and Nicolaus Copernicus (Toru\'n) for their support.



\footnotesize


\begin{thebibliography}{99}

\bibitem{BJL2016} Balan, R., Jakubowski, A. and Louhichi, S. (2016)
Functional convergence of linear processes with heavy-tailed innovations. J. Theoret. Probab. 29, 491--526.

\bibitem{bar2010} Barndorff-Nielsen, O. E. (2010) Stationary infinitely divisible processes. REBRAPE Braz. J. Probab. Stat.
25, 294--322.


\bibitem{bar2011} Barndorff-Nielsen, O. E., Benth, F. E. and Veraart, A. E. D.
 (2011) Recent advances in ambit stochastics. Preprint available at arXiv:1210.1354.


\bibitem{bar2014} Barndorff-Nielsen, O.E., Lunde, A., Shepard, N. and Veraart, A.E.D.
(2014) Integer-valued trawl processes:  a class of stationary infinitely divisible processes.
Scand.  J.  Statist. 41, 693--724.

\bibitem{Bill99} Billingsley, P.. (1999) {\em Convergence of Probability Measures}. 2nd ed., Wiley, New York.


\bibitem{fok} Christou, V.  and  Fokianos, K. (2014)
Quasi-likelihood inference for negative binomial time series.
J. Time Series Anal. 35, 55--78.



\bibitem{dav1970} Davydov, Y. A. (1970) The invariance principle for stationary processes.
Theor. Probab. Appl. 15, 487--498.




\bibitem{deh2002} Dehling, H. and Philipp, W. (2002) Empirical process techniques for dependent data. In:
H. Dehling, T. Mikosch and M. S{\o}rensen (Eds.),
{\em Empirical Process Techniques for Dependent Data}, pp. 1--113.
Birkh{\" a}user, Boston.












\bibitem{DOT03}
Doukhan, P. , Oppenheim, G. and  Taqqu  M. S. (Eds.)(2003)  {\em Theory and Applications of Long-Range Dependence}.
Birkh\"auser, Boston.

\bibitem{EPW1967} Esary, J.D., Proschan, F. and Walkup, D.W. (1967) Association of random variables, with applications. Ann. Math. Statist. 38, 1466--1474.


\bibitem{fell1966} Feller, W. (1966) {\em An Introduction to Probability Theory and Its Applications}, vol. 2. Wiley,
New York.











\bibitem{gir2012} Giraitis, L., Koul, H. L. and  Surgailis, D.   (2012) {\em Large Sample Inference
for Long Memory Processes.} Imperial College Press, London.

\bibitem{ibr1971} Ibragimov, I.A. and Linnink, Y.V. (1971)
{\em Independent and Stationary Sequences of Random Variables.} Wolters-Noordhoff, Groningen.

\bibitem{Jakubowski1997} Jakubowski, A. (1997)
A non-Skorohod topology on the Skorohod space. Electron. J. Probab. 2, 
1--21.


\bibitem{Hall1997} Hall, P., Koul, H.L. and Turlach, B.A. (1997)
Note on convergence rates of semiparametric estimators of dependence index.
Ann. Statist. 25, 1725--1739.



\bibitem{kajt2008} Kaj, I. and Taqqu, M. S. (2008)
Convergence to fractional Brownian motion and to the Telecom process: the integral representation approach.
In: Vares, M.E. and Sidoravicius, V. (Eds.) {\em An Out of Equilibrium 2.}
Progress in Probability, vol. 60, pp. 383--427. Birkh{\"a}user, Basel.


\bibitem{Konst1998}
Konstantopoulos, T. and Lin, S.-J. (1998) Macroscopic models for long-range dependent network traffic. Queueing Systems 28, 215--243.



\bibitem{kwa1992} Kwapie\'n, S. and Woyczy\'nski, W. A. (1992) {\em Random Series and
Stochastic Integrals: Single and Multiple.} Birkh\"auser, Boston.



\bibitem{ls2003} Leipus, R. and  Surgailis, D. (2003)  Random coefficient autoregression,
regime switching and long memory.
Adv. Appl. Probab. 35, 737--754.


\bibitem{lps2005} Leipus, R., Paulauskas, V. and  Surgailis, D. (2005)
Renewal regime switching and stable limit laws. J. Econometrics 129, 299-327.



\bibitem{lif2014} Lifshits, M. (2014) {\em Random Processes by Example.} World Scientific,  New Jersey.

\bibitem{LouhichiRio2011} Louhichi, S. and Rio, E. (2011) Functional convergence to stable L\'evy motions for iterated random Lipschitz mappings.
Electron. J. Probab. 16, 2452--2480.

\bibitem{miko2002} Mikosch, T., Resnick, S., Rootz\'en, H. and Stegeman, A. (2002)
Is network traffic approximated by stable L\'evy motion or fractional Brownian motion?
Ann. Appl. Probab. 12, 23--68.



\bibitem{pils2014} Pilipauskait\.e, V.  and Surgailis, D. (2014) Joint temporal and contemporaneous aggregation
of random-coefficient  AR(1) processes. Stochastic Process. Appl. 124, 1011--1035.



\bibitem{ResnickVdB2000} Resnick, S. and Van den Berg, E. (2000)  Weak convergence of high-speed traffic models.
J. Appl. Prob. 37, 375--397.

\bibitem{Skorokhod1956} Skorohod, A.V. (1956) Limit theorems for stochastic processes. Theory Probab. Appl. 1, 261--290.

\bibitem{su2004} Surgailis, D. (2004) Stable limits of sums of bounded functions of long memory
moving averages with finite variance.  Bernoulli 10, 327--355.

\bibitem{TaL1986}  Taqqu, M.S. and Levy, J.B. (1986)
Using renewal processes to  generate long-range dependence and high variability.
In: Eberlein, E. and Taqqu, M.S. (Eds.)
{\em Dependence in Probability and Statistics}, pp. 51--72.
Birkh\"auser, Boston.

\bibitem{TaWS1997}  Taqqu, M.S., Willinger, W. and Sherman, R. (1997) Proof of the
fundamental  result in self-similar traffic modeling. Computer Commun. Rev. 27, 5--23.

\bibitem{Wi2003} Willinger, W., Paxon, V., Riedi, R.H. and Taqqu, M.S. (2003)
Long-range dependence and data network traffic. In: Doukhan, P., Oppenheim, G. and Taqqu, M.S. (Eds.)
{\em Theory and Applications of Long-Range Dependence}, pp. 373--407.  Birkh\"auser, Boston.



\bibitem{wol2005} Wolpert, R. L. and Taqqu. M. S. (2005) Fractional Ornstein-Uhlenbeck L\'evy processes
and the Telecom process: upstairs and downstairs. Signal  Process. 85, 1523--1545.












\end{thebibliography}

\end{document}